\documentclass[letterpaper,12pt]{amsart}
\usepackage[backend=biber,maxcitenames=5,maxbibnames=9,doi=false,isbn=false,url=false,eprint=false, style=numeric]{biblatex}
\usepackage[margin=1.0in]{geometry}
\usepackage{graphicx}
\usepackage[shortlabels]{enumitem}
\usepackage[dvipsnames]{xcolor}
\usepackage{meyermacros}
\usepackage{ragged2e}

\renewbibmacro{in:}{}
\binoppenalty=\maxdimen
\relpenalty=\maxdimen

\title{Extending Fibrations of the $3$-Torus and Applications to Torus Surgery in $4$-Manifolds}
\author{Nicholas Meyer}
\address{Department of Mathematics and Statistics\\ \mbox{ Winona State University }\\ Winona, MN, United States}
\email{nick.meyer@winona.edu}
\urladdr{https://nickmeyer.phd}
\date{}
\subjclass{Primary: 57K40; Secondary: 57R22}
\keywords{fibration, torus surgery, knot surgery}

\usepackage{hyperref}
\hypersetup{
    pdfauthor={Nicholas Meyer},
    pdftitle={Extending Fibrations of the 3-Torus and Applications to Torus Surgery in 4-Manifolds}
}

\begin{document}


\begin{abstract}
Suppose that $W$ and $W'$ are smooth, compact, and oriented $4$-manifolds that are either diffeomorphic to $S^1$ times the exterior $E_Y(K)$ of a fibered knot $K$ in a closed, connected, orientable $3$-manifold $Y$, or are diffeomorphic to $\Sigma_{g,1}$ bundles over the $2$-torus with monodromy fixing the boundary of the fiber pointwise. If $f: \partial W' \to \partial W$ is an orientation-preserving diffeomorphism of the $3$-torus boundaries, we have that $X = W \cup_f W'$ is a closed, oriented $4$-manifold that fibers over $S^1$.
In particular, if $W' = T^2 \times D^2$ and $W = S^1 \times E_Y(K)$, then our result shows that the result of doing torus surgery in $S^1\times Y$ along $S^1 \times K$ is a $4$-manifold that fibers over $S^1$. Furthermore, we extend work of Zentner by showing that the result of torus surgery along $S^1$ times the unknot $\mathcal{U}$ in $S^1 \times S^3$ is diffeomorphic to $S^1$ times a lens space.
\end{abstract}

\maketitle


\section{Introduction}\label{sec:intro}

The study of $4$-manifolds via surgery techniques has a long and rich history, with connections to gauge theory, fibered knots, and surface bundles. One particularly fruitful method is \emph{torus surgery}, in which a $4$-manifold is modified by removing a tubular neighborhood of a smoothly embedded $2$-torus and re-gluing via a boundary diffeomorphism. This procedure generalizes Dehn surgery to dimension four and has been used to construct exotic smooth structures, as in the work of Dolgachev \autocite{Dolgachev2010} and Donaldson \autocite{Donaldson87}, wherein torus surgery produces infinite families of exotic elliptic surfaces $\CP^2 \cs 9 \overline{\CP}^2$.

Despite these advances, relatively little is known about how torus surgery interacts with fibrations of $4$-manifolds. Zentner \autocite{Zentner} showed that torus surgery along $S^1$ times the unknot, $\calU$, in $S^1 \times S^3$ preserves the diffeomorphism type if the resulting manifold is a homology $S^1 \times S^3$. This illustrates the rigidity of certain fibrations under surgery, but the general principles governing such behavior remain limited. In this paper, we extend these results and develop a broader framework for understanding how torus surgery preserves or modifies fiberedness.

A key ingredient in our approach is the use of \emph{fibered knot exteriors} in $3$-manifolds. Recall that a knot $K$ in a closed, connected, oriented $3$-manifold $Y$ is \emph{fibered} if its complement $E_Y(K) = Y \setminus \mathring{\nu}K$, where $\mathring{\nu}K$ denotes the interior of the closed tubular neighborhood of $K$ in $Y$, admits a smooth fibration over $S^1$, with fiber a Seifert surface $F$ and monodromy $\varphi$. The product $S^1 \times E_Y(K)$ then forms a natural $4$-manifold with $3$-torus boundary. By analyzing how fibrations extend from the boundary to the interior, we construct fibered $4$-manifolds arising from gluing knot exteriors to other $4$-manifolds, such as $T^2 \times D^2$ or other fibered manifolds with $T^3$ boundaries. This perspective generalizes classical \mbox{Fintushel-Stern} knot surgery \autocite{FS98} and allows for handling null-homologous knots in arbitrary $3$-manifolds while preserving fiberedness.

\mbox{Fintushel-Stern} knot surgery is a widely studied procedure for constructing exotic smooth structures on $4$-manifolds \autocites{Akbulut2002,AEMS2008,Yun2008,Ni17,Kim19}. It involves removing a tubular neighborhood of a self-intersection zero torus in a $4$-manifold $X$ and gluing in $S^1 \times E_{S^3}(K)$ via a diffeomorphism that preserves the homology class of the meridian. We adopt a slightly generalized version of this construction, allowing for null-homologous knots in arbitrary $3$-manifolds. We show that, under suitable conditions, the resulting $4$-manifold fibers over the circle. This unifies \mbox{Fintushel-Stern} knot surgery with torus surgery, highlighting the role of fibrations in controlling the topology of the resulting manifold.

In addition to products of a circle with fibered knot exteriors, we also consider \emph{surface bundles over the $2$-torus}. Let $\Sigma_{g,1}$ denote a compact, oriented surface of genus $g$ with one boundary component. A $\Sigma_{g,1}$ bundle over $T^2$ is a $4$-manifold $W$ equipped with a smooth submersion $W \onto T^2$ whose fibers are copies of $\Sigma_{g,1}$ and whose monodromy fixes the boundary pointwise. These manifolds have $3$-torus boundaries, and fibrations on the boundary extend across the interior using techniques similar to those for $S^1 \times E_Y(K)$. By treating both classes of $4$-manifolds --- $S^1$ times fibered knot exteriors and $\Sigma_{g,1}$ bundles over $T^2$ --- in a unified framework, we provide a systematic approach for constructing fibered $4$-manifolds arising from torus or generalized \mbox{Fintushel-Sten} knot surgery.

The central technique underlying our results is the extension of fibrations from the boundary $3$-torus to the interior of the manifold. We establish fibration extension theorems showing that if $W$ is either $S^1 \times E_Y(K)$ or a $\Sigma_{g,1}$ bundle over $T^2$, then many fibrations of $\partial W \cong T^3$ extend over $W$. These theorems are crucial for our analysis of torus surgery and generalized \mbox{Fintushel-Sten} knot surgery.


\section*{Acknowledgments}
This work was supported in part by the National Science Foundation under Grant No.\ DMS-2405301. The author would like to thank Rom\'{a}n Aranda and Alex Zupan for helpful conversations related to this work.


\section{Extending Fibrations}\label{sec:ExtFib}

In this section, we state and prove a series of theorems that allow us to extend a fibration of the $3$-torus across a $4$-manifold. These theorems will be used to construct fibrations that are compatible with various surgery constructions and will be used to construct new fibered $4$-manifolds.

The first case we investigate is the case of torus surgery along $S^1 \times K$ in $S^1 \times Y$, where $K$ is a fibered knot in a closed, connected, and oriented $3$-manifold $Y$. To do this, we establish a fibration extension theorem for $S^1$ times fibered knot exteriors.

\begin{theorem}\label{thm:K_ext}
	Let $Y$ be a closed, connected, oriented $3$-manifold, and let $K$ be a fibered knot in $Y$ with (oriented) fiber surface $F$ and monodromy $\varphi$. Define $W = S^1 \times E_Y(K)$ so that $\partial W = T^3$. Furthermore, let $\pi: T^3 \onto S^1$ be a fibration with $T$ being the fiber of $\pi$.

	If $T$ contains the curve $\lambda = \set{\mathrm{pt}} \times \partial F$, then $\pi$ extends to a fibration of $W = S^1 \times E_Y(K)$. That is, there is a fibration $\hat\pi: W \onto S^1$ such that $\restrict{\hat\pi}_{\partial W} = \pi$.

	\begin{proof}
		View $S^1 \times E_Y(K)$ as $S^1$ times an $F$-bundle over $S^1$, and let $S = S^1 \times \mu$ where $\mu$ is the meridian of $K$ that intersects the fiber surface $F$ once. The torus $S$ is an essential torus in $T^3$, for it is $\pi_1$-injective.

		We consider two cases, one where $S$ and $T$ are isotopic, and one where they are not. The case in which $S$ and $T$ are isotopic is less involved, so we will handle it first.

		If $S$ and $T$ are isotopic, then there is a curve $\alpha$ in $T^3$ that meets $S=T$ in a single point. Then we have a homology basis $\set{\lambda, \mu, \alpha}$ for $H_1(T^3)$.

		Now, assume that $S$ and $T$ are not isotopic. Up to isotopy, $T$ intersects $S$ in some family of simple closed curves, call these curves $\gamma_1, \dots, \gamma_n$. Observe that $\lambda$ intersects $S$ in a single point, as $\mu$ intersects $F$ in a single point. Now, since $T$ contains $\lambda$, we must have that $T$ intersects $S$ in a single curve, $\gamma$. Now, $\gamma$ must intersect $\lambda$ in a single point on $T$, so that $\lambda$ and $\gamma$ span $H_1(T)$.
		Now, pick a curve $\alpha$ in $S$ so that $\gamma$ and $\alpha$ intersect in a single point. This then gives that the set $\set{ \gamma, \lambda, \alpha}$ is a basis for $H_1(T^3)$.

		In either case, we have found a parameterization of $T^3$ via the homology bases we have found. In these coordinates, there is a fiber bundle map $\pi$  that is just projection onto the $\alpha$ coordinate, so that the pre-image of a point is a fiber isotopic to $T$.

		Now, we can see that there is a fiber bundle map $\hat\pi: S^1 \times E_Y(K) \onto S^1$, where we can parameterize $S^1 \times E_Y(K) = S^1_\gamma \times (F \tilde{\times}_\varphi S^1_\alpha)$ and take $\hat\pi$ to be projection onto the $S^1_\alpha$ factor. Then $\restrict{\hat\pi}_\partial = \pi$ as desired.
	\end{proof}
\end{theorem}

We next establish an analogous result for $\Sigma_{g,1}$ bundles over the $2$-torus. As before, the boundary of such a bundle is a $3$-torus, since the monodromy fixes the boundary of the fiber pointwise, and the key input is identifying the essential torus determined by a boundary section.

\begin{theorem}\label{thm:Surf_ext}
	Let $W$ be a $\Sigma_{g,1}$ bundle over $T^2$ so that $\partial W = T^3$, where the monodromy, $\rho$, of the fibration fixes the boundary $\partial \Sigma_{g,1}$ pointwise, and let $q: W \onto T^2$ be the bundle map. Furthermore, let $\pi: T^3 = \partial W \onto S^1$ be a fiber bundle with $T$ a fiber of $\pi$.

	If $T$ contains $\lambda = \partial \Sigma_{g,1} \times \set{\mathrm{pt}} $ then $\pi$ extends over $W$.

	\begin{proof}
		First, note that $\restrict{q}_{\partial W}: T^3 \to T^2$ is a fibration with $S^1$ fibers. Pick a section $s: T^3 \to T^2$ of this fibration, and let $S$ be the embedded image of $s$ under this section.

		The torus $S$ is an essential torus in $T^3$ since $q_*: \pi_1(T^3) \to \pi_1(T^2)$ is surjective and $q_*$ restricted to $S$ is an isomorphism.

		We work in two cases: one where $S$ and $T$ are isotopic, and one where they are not. The case in which $S$ and $T$ are isotopic is less involved, so we will handle it first.

		If $S$ and $T$ are isotopic, then there is a curve $\alpha$ in $T^3$ that meets $S$ or $T$ in a single point. Picking a homology basis $\gamma, \lambda$ for $H_1(S)\cong H_1(T)$,this basis extends to a homology basis $\gamma,\lambda, \alpha$ for $H_1(T^3)$.

		Now, if $S$ and $T$ are not isotopic, then up to isotopy, $T$ intersects $S$ in some family of simple closed curves, call them $\gamma_1, \dots, \gamma_n$. Observe that $\lambda $ intersects $S$ in a single point by construction. Now, since $T$ contains $\lambda$, we must have that $\abs{T \cap S} =1$. Let $\gamma$ be the curve $T \cap S$, and observe that $\gamma $ and $\lambda $ intersect in a single point on $T$ so that $\lambda$ and $\gamma$ span $H_1(T)$.
		Now, pick a curve $\alpha$ on $S$ so that $\alpha$ and $\gamma$ intersect in a single point.
		This then gives that the set $\set{ \gamma, \lambda, \alpha}$ is a basis for $H_1(T^3)$.

		Parameterize $T^3 = S^1_\gamma \times S^1_\lambda \times S^1_\alpha$. In these coordinates, the fiber bundle map $\pi$ is just projection onto the $\alpha$ coordinate, so that the pre-image of a point is a fiber isotopic to $T$.

		Now, we can see that there is a fiber bundle map $\hat\pi: W \to S^1$, where we can parameterize $W = \Sigma_{g,1} \tilde\times \pp{S^1_\gamma \times S^1_\alpha}$ and take $\hat\pi$ to be the projection onto the $S^1_\alpha$ factor. Then $\restrict{\hat{\pi}}_\partial = \pi$ as desired.
	\end{proof}
\end{theorem}

The preceding two results are central to this paper and will be used to construct many families of fibered $4$-manifolds.

The next lemma is a $3$-dimensional statement that underlies the applications of both of the above extension theorems. It characterizes essential $2$-tori in the $3$-torus from the point-of-view of fiber bundles and ensures that every essential torus arises as a fiber of a unique (up to diffeomorphism) trivial bundle.

\begin{lemma}\label{lem:essTorus}
	Let $T$ be an essential $2$-torus in $T^3$. Then there is a trivial fiber bundle $p: T^3 \to S^1$ with fibers isotopic to $T$, and the fiber bundle is unique up to diffeomorphism of $T^3$.
	\begin{proof}
		Pick a pair of dual curves $\alpha$, $\beta$ on $T$. Taking lifts to the universal cover, $\tilde\alpha$ and $\tilde\beta$ are isotopic to straight lines which span the plane $\tilde T$. This plane $\tilde T$ has a normal vector $\vec{n}$. The line spanned by $\vec{n}$ projects to a curve $\gamma$ in $T^3$ and $\gamma$ is dual to $T$. Parameterize $T^3$ as $S^1_\gamma \times T$, and let $p: T^3 \to S^1$ be projection onto the $S^1_\gamma$ factor. By construction, this is a fiber bundle and the fibers are isotopic to $T$. Moreover, this fiber bundle is trivial by construction.

		Now, suppose $p'$ is another fiber bundle of $T^3$ with fibers isotopic to $T$. In this case, pick a new basis $\alpha'$, $\beta'$ for a fiber $T'$ of $p'$. Since $T'$ is isotopic to $T$, we can assume that $\alpha'$ and $\beta'$ also live on $T$. In this case, there is a diffeomorphism $F: T^3 \to T^3$ sending the triple $(\alpha',\beta',\gamma')$ to $(\alpha, \beta,\gamma)$ since $\operatorname{Mod}(T^3) \cong \operatorname{GL}_3(\Z)$ and acts freely and transitively on the set of homology bases. In this way, we get that $p' = p \circ F$, as desired.
	\end{proof}
\end{lemma}

With these ingredients in place, we are now prepared to apply the extension theorems to construct large families of fibered $4$-manifolds arising from torus surgery and generalized \mbox{Fintushel-Sten} knot surgery.


\section{Gluing Theorems}\label{sec:GluingTheorems}

We first consider the case of gluing $S^1$ times fibered knot exteriors to $T^2 \times D^2$.

\begin{theorem}
	Let $K$ be a fibered knot in a closed, oriented $3$-manifold $Y$ with Seifert surface $F$ and monodromy $\varphi$. Let $W = T^2 \times D^2$, and let \mbox{$f: \partial W \to \partial(S^1 \times E_Y(K))$} be a diffeomorphism. Then, the manifold \mbox{$X_{K,f} = (S^1 \times E(K)) \cup_f W$} fibers over the circle.

	\begin{proof}
		Let $\lambda = \partial F \subset T^3 = \partial(S^1 \times E_Y(K))$, and let $\alpha = \set{\mathrm{pt}}  \times \partial D^2 \subset T^3 = \partial W$.
		If $\lambda$ and $f(\alpha)$ are not isotopic in $T^3$, then they span some essential $T^2$ in $T^3$.
		On the other hand, if $\lambda$ and $f(\alpha)$ are isotopic, there are infinitely many essential tori that contain $f(\alpha)$ and $\lambda$. In this case, let $T$ be one of them.

		In either case, by Lemma~\ref{lem:essTorus}, there is a fiber bundle $\pi: T^3 \to S^1$ whose fibers are isotopic to this essential torus.

		This fiber bundle map extends over $S^1 \times E(K)$ by Theorem~\ref{thm:K_ext}. The fiber bundle map $p$ also extends over \hbox{$W = T^2 \times D^2$} uniquely by Theorem~\ref{thm:K_ext} by observing that $T^2 \times D^2 = S^1 \times E_{S^3}(\calU)$.

		Therefore, $X_{K,f}$ is a fiber bundle over the circle.
	\end{proof}
\end{theorem}

We can also glue two fibered knot complements together to yield fibered manifolds:

\begin{theorem}
	Let $K_1, K_2$ be fibered knots in closed, oriented $3$-manifolds $Y_1, Y_2$, with Seifert surfaces $F_1, F_2$ and monodromies $\varphi_1, \varphi_2$. Let $W_i = S^1 \times E_{Y_i}(K_i)$ and let $f:\partial W_2 \to \partial W_1$ be a diffeomorphism (of $\partial W_1 = \partial W_2 = T^3$.)

	Let \mbox{$X_f = W_1 \cup_f W_2$.} Then $X_f$ fibers over the circle.

	\begin{proof}
		Let $\lambda_i = \partial F_i \subset T^3$.

		On one hand, if $\lambda_1$ and $f(\lambda_2)$ are not isotopic in $T^3$, then they span a unique essential $T^2$ in $T^3$ by lifting to the universal cover, straightening the curves, and noting that the plane spanned by the straightened lifts projects down to an essential torus.

		On the other hand, if $\lambda_1$ and $f(\lambda_2)$ are isotopic, there are infinitely many essential tori that contain $\lambda_1$ and $f(\lambda_2)$. In this case, let $T$ be one of them.

		In either case, by Lemma~\ref{lem:essTorus}, there is a fiber bundle $\pi: T^3 \to S^1$ whose fibers are isotopic to this essential torus.

		This fiber bundle map extends over each $W_i$ by Theorem~\ref{thm:K_ext}.
		Therefore, $X_{f}$ admits a fiber bundle over the circle.
	\end{proof}
\end{theorem}

As an immediate consequence, we now fully understand torus surgeries along tori of the form  $S^1 \times K$ in $S^1 \times Y$:
\begin{corollary}
	Let $K$ be a fibered knot in a closed, oriented $3$-manifold $Y$. Let $T = S^1 \times K \subset S^1 \times Y$. Then the result of doing torus surgery along $T$ in $X = S^1 \times Y$ fibers over the circle.
	\begin{proof}
		This is the previous theorem, applied to $T^2 \times D^2 = S^1 \times(S^1 \times D^2) = S^1 \times E_{S^3}(\calU)$.
	\end{proof}
\end{corollary}

We now extend Zentner's result\ \autocite{Zentner} on torus surgeries along $S^1\times\mathrm{U}$ in $S^1 \times S^3$:

\begin{corollary}\label{cor:lens_space}
	The result of torus surgery along $S^1\times \calU'$ in $S^1 \times S^3$ is diffeomorphic to $S^1$ times a lens space.
	\begin{proof}
		Let $X= (S^1\times S^3)_{S^1 \times \calU',f}$ be the surgery manifold. View $X = (T^2 \times D^2) \cup_f (T^2\times D^2)$. Let $\alpha = \set{\mathrm{pt}}  \times \partial D^2$ in the first factor, and let $\beta = \set{\mathrm{pt}}  \times \partial D^2 $ in the second. There is an essential $2$-torus $T$ in $T^3$ containing $\alpha$ and $f(\beta)$. Use Lemma~\ref{lem:essTorus} to construct a fibration $\pi: T^3 \to S^1$ whose fibers are isotopic to $T$ and the fibration map is projection onto some curve $\gamma$ in $T^3$. This fibration extends trivially over each factor $T^2 \times D^2$, and the fibers of $\hat\pi$ have a Heegaard splitting given by $H_\alpha = S^1 \times D^2_\alpha$ and $H_\beta = S^1\times D^2_\beta$, and the Heegaard surface is one of the parallel copies of $T$. These fibers are lens spaces, since this is a Heegaard splitting of genus one. Hence, $X$ is diffeomorphic to $S^1 \times Y$ where $Y$ is a lens space.
	\end{proof}
\end{corollary}

We now turn our attention to the surface bundle case, wherein we gently modify the proof from the knot complement case.

\begin{theorem}
	Let $W $ be a $\Sigma_{g,1}$ bundle over $T^2$ and let $W'$ be a $\Sigma_{g',1}$ bundle over $T^2$.  Suppose $f: \partial W' \to \partial W$ is a diffeomorphism of the $3$-torus, and set $X= W \cup_f W'$. Then $X$ fibers over $S^1$.

	\begin{proof}
		Let $\alpha = \partial \Sigma_{g,1} \subset T^3$, and let $\beta = \partial \Sigma_{g',1} \subset T^3$. If $\alpha$ and $f(\beta)$ are not isotopic in $T^3$, then they span a unique (up to isotopy) essential $2$-torus in $T^3$.
		On the other hand, if $\alpha$ and $f(\beta)$ are isotopic, there are infinitely many essential tori that contain $\alpha$ and $f(\beta)$. In this case, let $T$ be one of them.
		In either case, by Lemma~\ref{lem:essTorus}, there is a fiber bundle $\pi: T^3 \to S^1$ whose fibers are isotopic to this essential torus.
		This fiber bundle map extends over each $W_i$ by Theorem~\ref{thm:Surf_ext}.
		Therefore, $X$ admits a fibration over the circle.
	\end{proof}
\end{theorem}

Finally, we combine the knot-complement and surface-bundle cases. The proof is identical to the previous two theorems, with the essential torus argument being the key step yet again.

\begin{corollary}\label{cor:genExt}
	Suppose that $W$ and $W'$ are smooth, compact, and oriented $4$-manifolds that are diffeomorphic to $S^1 \times E_Y(K)$, where $Y$ is a closed, oriented $3$-manifold and $K$ is a fibered knot in $Y$, or to $\Sigma_{g,1}$ bundles over the $2$-torus.

	If $f: \partial W' \to \partial W$ is an orientation-preserving diffeomorphism of the $3$-torus boundaries, we have that $X = W \cup_f W'$ is a closed, oriented $4$-manifold that fibers over $S^1$.
\end{corollary}


\section{Examples}\label{sec:Examples}
We now give some examples of the preceding results. Our first example demonstrates that we can obtain $S^1$ times a lens space by performing torus surgery along $S^1 \times \calU$ in $S^1 \times S^3$.

\begin{example}\label{ex:lens_space}
	Let $T = S^1 \times \calU$ in $S^1 \times S^3$. Let $\lambda$ and $\mu$ be a longitude and meridian of $\calU$ in $S^3$, respectively. Then $\lambda$ and $\mu$ live as curves in $T^3 = \partial((S^1 \times S^3 )\setminus \mathring\nu T)$. The complement of $\nu T$ in $S^1 \times S^3$ is diffeomorphic to $D^2 \times T^2$. We can parameterize $\partial((S^1 \times S^3 )\setminus \mathring\nu T)$ as $D^2 \times S^1_\lambda \times S^1_\mu$.

	Let $f: \partial(D^2 \times T^2) \to \partial ((S^1 \times S^3) \setminus \mathring\nu T)$ be a diffeomorphism of $T^3$ that sends $\partial D^2$ to $p \lambda + q \mu$, and identifies one of the remaining $S^1$ factors of $\partial(D^2 \times T^2)$ with $S^1_\lambda$. Thus, the results of doing $f$-torus surgery along $T$ is diffeomorphic to the lens space $S^1 \times L(q,p)$.
\end{example}

Now, using Corollary \ref{cor:lens_space}, we obtain a complete diffeomorphism classification of all possible torus surgeries along $S^1 \times \calU$ in $S^1 \times S^3$.

We now introduce a framework for generalized \mbox{Fintushel-Stern} knot surgery, which replaces a torus neighborhood in a $4$-manifold with $S^1$ times a knot complement:

\begin{definition}[Generalized \mbox{Fintushel-Stern} Knot Surgery, adapted from \autocite{FS98}]
	Let $X$ be a closed, connected, oriented $4$-manifold, let $T$ be an embedded $2$-torus with self-intersection number zero, and let $K$ be a null-homologous knot in a closed $3$-manifold $Y$.

	Suppose that \mbox{$\varphi: \partial(\nu T)\to \partial(S^1 \times E_Y(K))$} sends $\set{\mathrm{pt}}\times \partial D^2$ to a longitude of $K$ in $S^1 \times E_Y(K)$. The $4$-manifold $X_{T, K, \varphi} = (X \setminus \mathring{\nu}T ) \cup_\varphi (S^1 \times E_Y(K))$ is said to be obtained by \emph{generalized \hbox{Fintushel-Stern} knot surgery along $T$ with knot $K$.}
\end{definition}

For the unknot $K = \calU$, this construction recovers the identity surgery: $X_{T,\calU,\varphi}$ is diffeomorphic to $X$ for any $T$ and $\varphi$ satisfying the definition above.

We focus next on a special class of tori that allow us to better understand some generalized \mbox{Fintushel-Stern} knot surgeries:

\begin{definition}\label{def:tft}
	Let $X$ be a smooth, closed $4$-manifold and let $T: T^2 \into X$ be an embedded $2$-torus. We say that $T$ is a \emph{torus-fibered $T^2$-knot} if $T$ has self-intersection number zero and the complement of the normal bundle of $T$ in $X$ fibers over the $2$-torus such that the monodromy of the fibration fixes the boundary of the fiber pointwise.
\end{definition}

Natural examples of torus-fibered $T^2$-knots arise from fibered knots in $3$-manifolds:

\begin{example}\label{ex:tft}
	Let $K$ be a fibered knot in a closed, connected, oriented $3$-manifold $Y$ with Seifert surface $F$ and monodromy $\varphi$. Then the 2-torus $T = S^1 \times K$ in $S^1 \times Y$ is a $T^2$-knot whose complement is fibered over the torus. To see this, take $p: E_Y(K) \to S^1$ to be a fibration on the complement of $K$ in $Y$, and define $\pi: S^1 \times E_Y(K) \to S^1 \times S^1 = T^2$ by $\pi(\theta, x) = (\theta, p(x))$. This map $\pi$ is a fibration over the torus, as desired.
\end{example}

This observation leads to a description of certain generalized \mbox{Fintushel-Stern} knot surgeries:

\begin{corollary}
	Let $T$ be a torus-fibered $T^2$-knot in a closed $4$-manifold $X$, and let $K$ be a fibered knot in a closed $3$-manifold $Y$. Suppose that $f: \partial(S^1 \times E_Y(K))\to \partial(X \setminus \mathring{\nu}T)$ is a diffeomorphism of the $3$-torus that sends $\set{\mathrm{pt}} \times \lambda$ to $\set{\mathrm{pt}} \times \partial D^2$, where $\lambda$ is a longitude of $K$ in $E_Y(K)$. Then the generalized \mbox{Fintushel-Stern} knot surgery $X_{T,K,f}$ along $T$ with knot $K$ and mapping $f$ fibers over the circle.
	\begin{proof}
		Apply Corollary \ref{cor:genExt} to the case where $W$ is $S^1 \times E_Y(K)$ and $W'$ is $X \setminus \mathring{\nu}T$, which is a $\Sigma_{g,1}$ bundle over the $2$-torus, since $T$ is a torus-fibered $T^2$-knot.
	\end{proof}
\end{corollary}

An algebro-topological constraint for the existence of torus-fibered $T^2$-knots is as follows:

\begin{theorem}\label{thm:t2fibalg}
	If $X$ admits a torus-fibered $T^2$-knot, then the Euler characteristic $\chi(X)$ and signature $\sigma(X)$ both vanish.
	\begin{proof}
		Using the inclusion-exclusion theorem for Euler characteristic and the fact that the Euler characteristic of a bundle is multiplicative, we have the following:
		\begin{align*}
		\chi(X) &= \chi(T^2 \times D^2) + \chi(E_{X}(T))- \chi(T^3) \\
		&= 0 + \chi(\Sigma_{g,1}\tilde\times_\rho T^2) - 0 \\
		&= \chi(\Sigma_{g,1}) \cdot \chi(T^2)\\
		&= 0.
		\end{align*}

		We now compute the signature result. Since $E_X(T)$ fibers over the 2-torus, $E_X(T)$ also fibers over $S^1$: pick non-null-homotopic curve $\gamma$ on the base space $T^2$, and compose the fibration map with the projection onto $\gamma$. This shows that $X$ admits an open book decomposition with binding the torus-fibered $T^2$-knot, so its signature must vanish (see~\autocite{Ranicki1998}, Chapter 29).
	\end{proof}
\end{theorem}

This theorem shows that torus-fibered $T^2$-knots impose strong topological restrictions on the ambient $4$-manifold. While necessary, these conditions are not always sufficient: an open book decomposition with torus binding does not guarantee the existence of a torus-fibered $T^2$-knot, since the page may fail to be a surface bundle or the monodromy may not fix the boundary pointwise.

Overall, these examples unify the theme of the section: torus surgery, particularly torus surgery on torus-fibered $T^2$-knots, provides a concrete, bundle-theoretic pathway to construct and understand generalized \mbox{Fintushel-Stern} knot surgeries in $4$-manifolds.

In light of Theorem	\ref{thm:t2fibalg}, it is natural to ask which $4$-manifolds with vanishing Euler characteristic and signature admit torus-fibered $T^2$-knots. In particular, do $\cs^k (S^1 \times S^3)$ admit torus-fibered $T^2$-knots for all $k\geq 2$? (Example \ref{ex:tft} is the case where $k=1$.) We leave this question for future work. Moreover, in light of the proof of Theorem \ref{thm:t2fibalg}, it is natural to ask whether every $4$-manifold that admits an open book decomposition with torus binding also admits a torus-fibered $T^2$-knot. We also leave this question for future work.


\clearpage
\printbibliography

\end{document}